\documentclass[12pt]{iopart}
\pdfoutput=1
\usepackage{graphicx, amsfonts, amsopn, amsbsy, 
            amssymb, amsthm, verbatim, latexsym}

\newcommand{\dd}{{\rm d}}
\newcommand{\id}{{\rm id.}}

\newcommand{\one}{\underline{1}}

\newcommand{\bC}{{\mathbb{C}}}

\newcommand{\bR}{{\mathbb{R}}}

\newcommand{\bN}{{\mathbb{N}}}

\newcommand{\rR}{{\mathrm{R}}}

\newcommand{\mean}[1]{\overline{#1}}

\newcommand{\ave}[1]{\left\langle #1 \right\rangle}

\newcommand{\beq}{\begin{equation}}
\newcommand{\eeq}{\end{equation}}
\newcommand{\beqa}{\begin{eqnarray}}
\newcommand{\eeqa}{\end{eqnarray}}
\newcommand{\beqas}{\begin{eqnarray*}}
\newcommand{\eeqas}{\end{eqnarray*}}

\begin{document}

\title{The ensemble of random Markov matrices}

\author{Martin Horvat}

\address{Physics Department, Faculty of Mathematics and Physics, University of Ljubljana, Slovenia}

\address{CAMTP - Center for Applied Mathematics and Theoretical Physics, University of Maribor, Krekova 2, SI-2000 Maribor, Slovenia}

\eads{\mailto{martin.horvat@fmf.uni-lj.si}}
\begin{abstract}
The ensemble of random Markov matrices is introduced as a set of Markov or stochastic matrices with the maximal Shannon entropy. The statistical properties of the stationary distribution $\pi$, the average entropy growth rate $h$ and the second largest eigenvalue $\nu$ across the ensemble are studied. It is shown and heuristically proven that the entropy growth-rate and second largest eigenvalue of Markov matrices scale in average with dimension of matrices $d$ as $h\sim \log(O(d))$ and $|\nu| \sim d^{-\frac{1}{2}}$, respectively, yielding the asymptotic relation $h \tau_{\rm c}\sim 1/2$ between entropy $h$ and correlation decay time $\tau_{\rm c} =-1/\log|\nu|$ . Additionally, the correlation between $h$ and and $\tau_{\rm c}$ is analysed and is decreasing with increasing dimension $d$.
\end{abstract}

\submitto{{\it J. Statistical Mechanics}}

\section{Introduction}
In the information theory and mathematical modeling of physical processes we often stumble upon Markov chains \cite{meyn:book:95} and Markov or stochastic matrices \cite{meyer:book:00}, which determine the evolution of the first. Let us assume that a Markov chain is based on the set of states ${\cal S} =\{s_i\}_{i=1}^d$ and  $P(s_j|s_i)$ is the conditional probability for a transition from state $s_i$ to $s_j$ then the corresponding Markov matrix $M\in\bR_+^{d\times d}$ is a collection of conditional probabilities
$$
  M_{i,j} = P(s_j|s_i)\ge 0 \>, \qquad \sum_{j=1}^d M_{i,j} = 1\>,
$$
where $d$ the dimension of the Markov matrix.  Notice that the sum of elements in each rows is normalized to 1. The applications of Markov matrices as well as the ways how to constuct them are very diverse. In particular interesting is their use in dynamical systems for a probabilistic description of the dynamics. For a general introduction in this direction see e.g. \cite{christiansen:non:96}. For example let us consider a discrete dynamical system $\phi^t:X\to X$, where $t\in\bN$, with the phase space $X$ and the invariant measure $\mu$. By choosing disjoint subsets of phase space $\{X_i \subset X:X_i\cap X_j=0~ \textrm{for}~ i\neq j\}$, which satisfy $\bigcup_i X_i =X$, the Markov matrix $M=[M_{i,j}]_{i,j=1}^d$ corresponding to the dynamical system can be defined as 
$$
  M_{i,j} = \frac{\mu(\phi(X_i) \cap X_j)}{\mu(X_i)}\>,
$$
and describes a single time step of the dynamical system. In this way a paradigmatic example of a dynamical system with an algebraic decay of correlation -- the triangle map \cite{casati:prl:00} was examined in reference \cite{horvat:physD:09}. Beside the presented method to construct a Markov matrix and other methods producing matrices with specific properties, there is often a need to construct a Markov matrix ad-hoc i.e. without incorporating any information about the system except the number of states $d$. If the construction procedure is a stochastic process then the resulting matrix is called {\it the random Markov matrix} and the set of such matrices form {\it the ensemble of random Markov matrices}. These matrices are usually used, without much theoretical background, for testing purposes. For example testing of algorithms or certain statistical hypothesises, where application to the field of dynamical systems, connected to ergodicity and mixing property, are the most interesting to us. In the information theory random Markov matrices are used to test the algorithms of recognition or attribution process, compression algorithms etc.\par
The presented work is strongly related to the work of \cite{chafai:07} discussing the general properties of the Markov ensemble and to references \cite{goldberg:meth:00} and \cite{goldberg:siam:03}, where a closer look at the second largest (sub-dominant) eigenvalue of random Markov matrices was taken. In contrast to the past works we try to emphasise the physical application of results in particular in the field of dynamical systems. Interestingly, recently ensemble of random Markov matrices was applied in the spectral analysis of the random quantum super-operators \cite{zyczkowski:preprint}.

\section{Preliminaries}

The set of all possible Markov matrices $M$ of dimension $d$ is defined as
$$
  {\cal M}(d) = \{ M \in \bR_+^{d\times d}: M \one = \one \}\>,\qquad
  \one = (1,\ldots, 1)\in \bR^d\>,
$$
and it is isomorphic to the direct product of $d$ convex sets
$$
 \bigotimes^d \{ x\in \bR_+^d ~:~ x^T \one = 1\} \>.
$$
The set ${\cal M}(d)$ forms together with the matrix product a semi-group, whereas the set of non-singular Markov matrices form a group of stochastic matrices. The ensemble of random Markov matrices is defined as a set of Markov matrices ${\cal M}(d)$ with the probability measure of matrix elements $M_{i,j}\in\bR_+$ reading
\beq
  \dd P(M) = [(d-1)!]^d \delta^d(M\one - \one) \dd M\>,\qquad 
  \dd M := \prod_{i,j=1}^d \dd M_{i,j} \>,
  \label{eq:distr_mat}
\eeq
which incorporates minimal information about the set i.e. only constrains due to the probability conservation in the Markov process. The ensemble of random Markov matrices, denoted by the pair $(\dd P(M), {\cal M}(d))$, is also referred to as the Dirichlet ensemble \cite{chafai:07} and corresponding matrices are called doubly (row and column) stochastic or bi-stochastic Markov matrices \cite{meyer:book:00}. \par
The rows of the Markov matrix from the ensemble $(\dd P(M), {\cal M}(d))$ are independent random vectors $X=(X_i\ge 0)_{i=1}^d$ with the distribution
\beq
  P_{\rm rows} (X) = (d-1)!\, \delta(\one^T X-1)\>.
 \label{eq:P_rows}
\eeq
It can be rather awkward to numerically generate components of vector-rows $X$ directly and so a different approach to do that is taken. By following the notes on the exponential distribution in \cite{feller:book:66}, p. 76, we find that the vector-rows $X$ of Markov matrices can be expressed by vectors $Y=(Y_i)_{i=1}^d$ of independent variables $Y_i$ with a common exponential distribution in the following way
\beq
 X = \frac{Y}{\one^T Y}\>.
 \label{eq:gen_rows}
\eeq
This the way we numerically generate the pseudo-random rows in a Markov matrices, where each row generated independently. Consequently, the distribution of the rows can be written as
\beq
  P_{\rm rows} (X)  =
  \int_{\bR_+^d} \dd^d Y \delta^d\left(X - \frac{Y}{\one^T Y}\right) e^{-\one^T Y} \>.
  \label{eq:P_rows_1} 
\eeq
This identity of expressions is checked by calculating the moment-generating function the expression above and expression (\ref{eq:P_rows}) yielding the same result equal to
$$
  \int_{\bR_+^d} \dd^d X e^{-\lambda^T X} P_{\rm rows}(X)
  = (-1)^{d-1}(d-1)!\,  \sum_{i=1}^d \frac{e^{-\lambda_i}}{w'(\lambda_i)} \>,
$$
where $\lambda=(\lambda_i)_{i=1}^d$ and $w(x)= \prod_{i=1}^d(x-\lambda_i)$. Let us denote the sum of variables $Y_i$ by $S=\one^T Y$ and examine its statistics. The ratio between the standard deviation $\sigma_S=\sqrt{\ave{(S - \ave{S})^2}}$ of the $S$ and its average value $\ave{S}$ is equal to $\sigma_S/\ave{S}=d^{-\frac{1}{2}}$ and it is decreasing with increasing dimension $d$.  This means that the renormalization of variables $Y_i$ in the expression (\ref{eq:gen_rows}) has less and less effect on the functional dependence between $X$ and $Y$, as $d$ is increased. We conclude that in the limit for large dimensions $d\gg 1$ variables $X_i$ are approximately independent and exponentially distributed with distribution $P_{\rm m}(X_i)=d\, \exp(-d\, X_i)$. Following this idea we write the asymptotic approximation of the probability measure of Markov ensemble as
\beq
  \dd P(M) \sim 
  \dd P_{\rm asym}(M) = d^{d^2}e^{-d\, \one^T M \one}\, \dd M \>,
  \label{eq:distr_mat_asym}
\eeq
which is valid in the limit $d\to\infty$. It can be verified that the averages w.r.t. the distributions $\dd P(M)$ and $\dd P_{\rm asym}(M)$ of an well behaved observable on $\cal M$, which depends only on a finite number of matrix elements, are asymptotically equal.\par
The probability measure of the Markov ensemble $\dd P(M)$ has several unique properties that makes the defined ensemble of random Markov matrices $(\dd P(M), {\cal M}(d)$ interesting and potentially useful.\par
For instance the probability measure $\dd P(M)$ has a maximal Shannon entropy and in this information sense it is unique. The set of Markov matrices is merely a direct product of planes restricted to $\bR_+^d$ and $\dd P(M)$ is uniform on them. The Shannon entropy of the measure $\dd P(M)$ on the set ${\cal M}(d)$ is just the sum of Shannon entropies of the uniform distribution on the planes, which are themselves minimal. Hence the Shannon entropy of $\dd P(M)$ is also minimal. Any modification of the measure would necessarily increase the Shannon entropy and therefore it is unique.\par
It is also interesting to notice that the measure $\dd P(M)$ is not invariant w.r.t. matrix multiplication. However, for a given non-singular Markov matrix $A\in {\cal M}(d)$ the measure $\dd P(M)$ is invariant on matrix multiplication up a constant 
$$
   P(A {\cal B}) 
  = |\det (A)|^{-d}  P({\cal B})\qquad \forall {\cal B} \subset  {\cal M}(d)\>.
$$
In fact there is no measure of Markov matrices with the matrix elements $M_{i,j}$ approximately independent in the limit $d\gg 1$, which would be invariant w.r.t. matrix multiplication. To show this let us consider two large matrices $A=[A_{i,j}]_{i,j=1}^d$ and $B=[B_{i,j}]_{i,j=1}^d$ with the matrix elements being i.i.d. variables with the distribution $P_{\rm m}(x)$. Here we denote $i$-th central moment of some distribution $Q(x)$ for $i>1$ as $\mu_i(Q) = \int \dd x\,(x-\mu_1(Q))^i Q(x) $ and for $i=1$ by $\mu_1(Q)=\int \dd x\, Q(x) x$. We assume that the first three central moments of $P_{\rm m}(x)$ are finite. The matrix elements of the product $AB=[C_{i,j}=\sum_{k=1}^d A_{i,k} B_{k,j}]_{i,j=1}^d$ are distributed by
$$
  P_{\rm C}(x) = \underbrace{P_{\rm AB} * \ldots * P_{\rm AB}}_{d}(x)\>, \qquad 
  P_{\rm AB} (x) = \int_{\bR_+} \frac{\dd a}{a}\,P_{\rm m}(a)\, P_{\rm m}\left(\frac{x}{a}\right) \dd a\>, 
$$
where the sign $*$ denotes the convolution and $P_{\rm AB}$ is the distribution of the product $A_{i,j} C_{j,k}$ with the first two central moments reading: 
$$
\mu_1(P_{\rm AB})=\mu_1(P_{\rm m})^2\quad\textrm{and}\quad
\mu_2(P_{\rm AB}) = 2\left(\mu_1(P_{\rm m})^2 + \mu_2(P_{\rm m})\right) \mu_2(P_{\rm m}) \>.
$$
According to the central limit theorem (CLT) \cite{feller:book:66} the distribution $P_{\rm C}(x)$ converges in the limit $d\to\infty$ to a Gaussian distribution with the first two central moments equal to $\mu_1(P_{\rm C})=d\,\mu_1(P_{\rm AB})$ and $\mu_2(P_{\rm C})=d\, \mu_2(P_{\rm AB})$. For distribution $P_{\rm m}(x)$ to be invariant w.r.t. to the matrix multiplication it has to be asymptotically equal to $P_{\rm C}(x)$ meaning that $P_{\rm m}(x)$ is also a Gaussian distribution for large dimension $d\gg 1$. By comparing the first two central moments of $P_{\rm C}(x)$ and $P_{\rm m}(x)$ we find that the average value of matrix elements of the Markov matrix and their variance are asymptotically equivalent to 
$$
  \mu_1(P_{\rm m})\sim \frac{1}{d}\quad\textrm{and}\quad
  \mu_2(P_{\rm m})\sim \frac{1}{d}-\frac{2}{d^2}\>,
$$
respectively. The ratio between the standard deviation and the average scales with dimension as $\sqrt{\mu_2(P_{\rm m})}/\mu_1(P_{\rm m})=O(d^{\frac{1}{2}})$ and diverges in the limit $d\to\infty$. This indicates that a measure of Markov matrices by which the matrix elements are asymptotically independent and distributed by $P_{\rm m}(x)$ does not exist. \par
In the following we discuss the properties of random Markov matrices from the Markov ensemble $(\dd P(M), \cal M)$. We focus on the entropy growth rate and correlation decay in the Markov chains generated by these Markov matrices, and examine their asymptotic behaviour for $d\gg 1$.

\section{The entropy and stationary distribution of the random Markov matrices}

\begin{figure}[!htb]
  \centering
  \includegraphics[width=7cm]{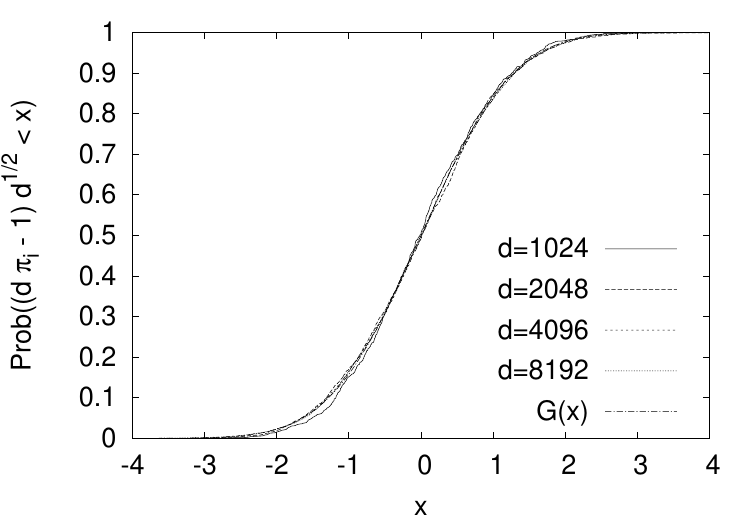}
  \includegraphics[width=7cm]{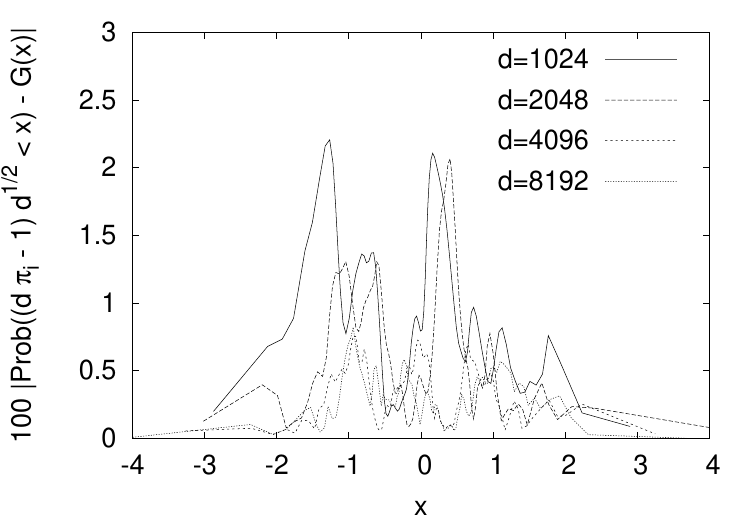}
\hbox to15cm{\scriptsize\hfil(a)\hfil(b)\hfil}
  \caption{The cumulative distribution of the rescaled coefficients $\pi_i$ corresponding to a stationary distributions $\pi=(\pi_i)_{i=1}^d$ (a) and its deviation form the expected limiting form $G(x)$ (\ref{eq:G_def}) (b) calculated 
for individual random Markov matrices of different dimensions $d$ taken from the ensemble.}
  \label{pic:inv_cum_dist}
\end{figure}

\begin{figure}[!htb]
  \centering
  \includegraphics[width=7cm]{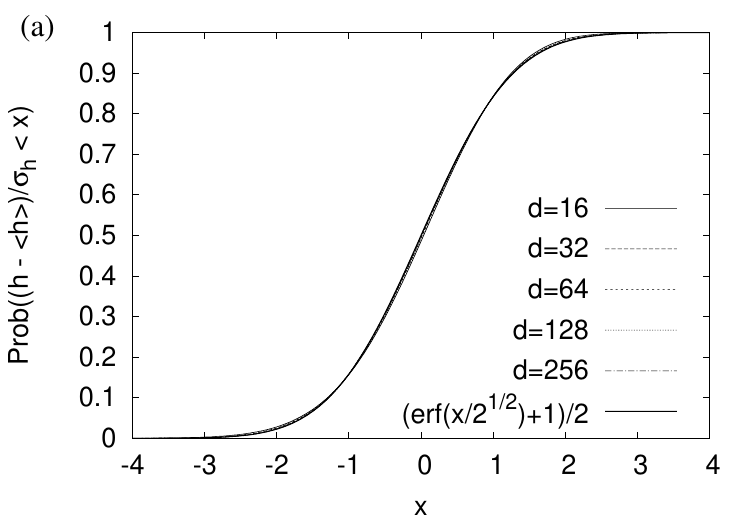}\\
  \includegraphics[width=7cm]{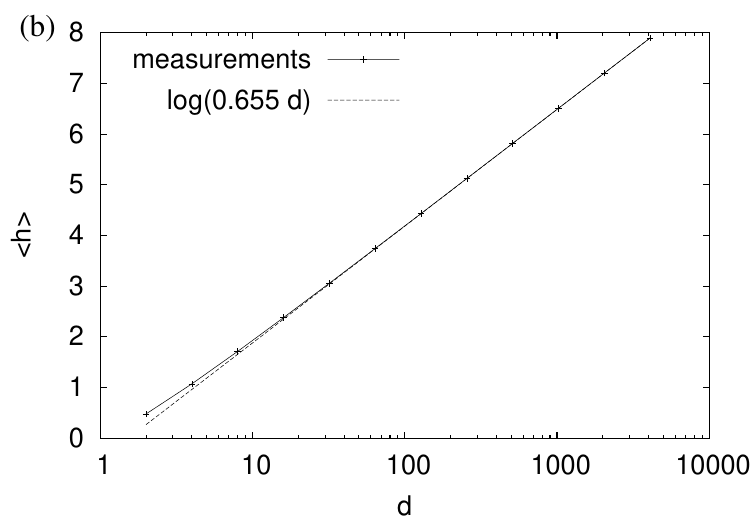}%
  \includegraphics[width=7cm]{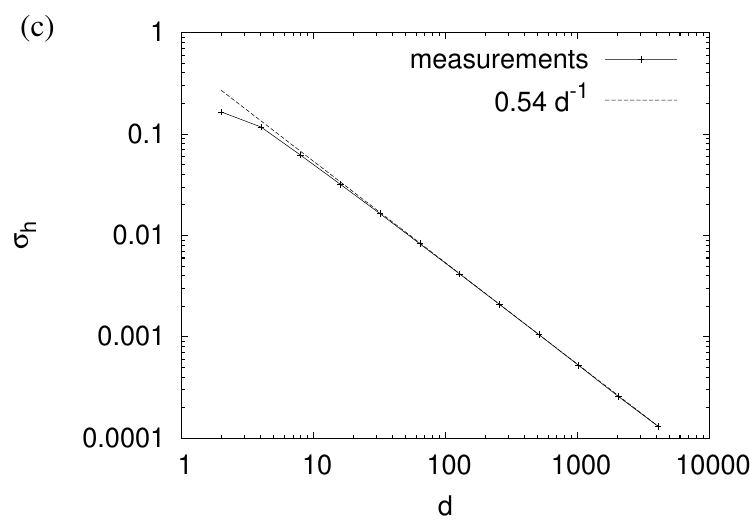}
  \caption{The cumulative distribution of the entropy growth rate $h$ of random Markov matrices for different dimensions $d$ (a), and the average entropy growth rate $\ave{h}$ (b) and its standard deviation $\sigma_h$ (c) as a function of dimension $d$ calculated using $N=10^6$ random Markov matrices from the ensemble.}
  \label{pic:entropy}
\end{figure}

We consider a Markov chain defined on the set of states ${\cal S}= \{s_i\}_{i=1}^d$ and with the conditional probabilities $P(s_j|s_i)=M_{i,j}$ given in the Markov matrix $M=[M_{i,j}]_{i,j=1}^d$. The initial probability distribution over the states is $(P(s_i))_{i=1}^d$. The probability that the Markov chain has evolved up to time $t$ following a specific route $(e_1,\ldots,e_t) \in {\cal S}^t$ is given with the product of conditional probabilities reading
$$
  P(e_1,\ldots, e_t) = P(e_1) P(e_2|e_1) P(e_3| e_2) \ldots P(e_t | e_{t-1})\>.
$$
Then the dynamic entropy $S$ of the Markov chain at time $t$ is given by the sum 
$$
  S(t) =- \sum_{e\in {\cal S}^t}  P(e)\log P(e)\>,
$$
taken over all different routes up to time $t$. In ergodic Markov chains we expect that the entropy in the limit $t\to\infty$ increases linearly with increasing time $t$ as
$$
  S(t)\sim h\, t \>,
$$ 
where $h\in \bR$ denotes the asymptotic entropy growth rate of the Markov chain. The entropy growth rate is given by the formula \cite{slomczy:open:02}
\beq
  h = -\sum_{i=1}^d \pi_i \sum_{j=1}^d M_{i,j} \log M_{i,j}\>,
  \label{eq:h_def}
\eeq
where we use the stationary distribution $\pi=(\pi_i\ge 0)_{i=1}^d$ defined as the eigenvector of the Markov matrix corresponding to the unit eigenvalue,
\beq
  \pi^T M = \pi^T\>, \qquad  \sum_{i=1}^d \pi_i =1 \>.
  \label{eq:inv}
\eeq
In the following we discuss the distribution $P_\pi(x)$ of elements $\pi_i$ corresponding to a stationary distribution $\pi$ of a random Markov matrix $M=[M_{i,j}]_{i,j=1}^d$. 
In particular we are interested in the asymptotic limit as $d\to\infty$, where the matrix elements $M_{i,j}$ are approximately independent variables with an exponential distribution $P_{\rm m}(M_{i,j})=d \exp(-d\,M_{i,j})$. Further we assume that $p_i$ and $M_{i,j}$ have no correlations. By doing this the eigenvalue equation (\ref{eq:inv}) written by components $\pi_i= \sum_j \pi_j M_{j,i}$ can be translated into an asymptotic self-consistent condition for the distribution $P_\pi(x)$ reading
\beq
  P_\pi(x) \sim \underbrace{P_{\rm \pi M} * \ldots * P_{\rm \pi M}}_{d}(x)
  \quad \textrm{as}\quad d\to\infty \>,
  \label{eq:pi_self}
\eeq
with the distribution $P_{\rm \pi M}(x)$ of the products $p_j M_{j,i}$ depending again on distribution $P_\pi(x)$ and written as
$$
  P_{\rm \pi M}(x) = d \int_{\rR_+} \frac{\dd b}{b}  
     \exp\left(-d\frac{x}{b}\right) P_\pi(b)\>,
$$
where $*$ denotes the convolution. Assuming that the first three central moments of the distribution $P_{\rm \pi M}(x)$ are finite we can use the CLT and state that the distribution $P_\pi(x)$ converges to a Gaussian distribution as $d\to\infty$. By inserting the ansatz 
$$
 P_\pi (x) = \frac{1}{\sqrt{2\pi \sigma_\pi^2}}
\exp\left(-\frac{(x-\mean{\pi})^2}{2\sigma_\pi^2}\right)\>,\qquad 
  \mean{\pi}=\frac{1}{d}\>,
$$
into equation (\ref{eq:pi_self}) and imposing the asymptotic equality (\ref{eq:pi_self}) we obtain the variance $\sigma_\pi^2\sim d^{-3}$ of the elements $\pi_i$. By appropriately  rescaling the coefficients $\pi_i$ their cumulative distribution is independent of dimension $d$ in the limit $d\to\infty$ and reads
\beq
  {\rm Prob} \left( \frac{d \pi_i - 1}{\sqrt{2}} <\frac{x}{\sqrt{d}}\right) 
  \sim G(x) = \frac{1}{2} \left({\rm erf}\left(\frac{x}{\sqrt{2}}\right) + 1\right)\>.
  \label{eq:G_def}
\eeq
This result is compared in figure \ref{pic:inv_cum_dist} with numerically obtained distributions of rescaled coefficients $\pi_i$ for different large dimensions $d$ and we find a very good agreement.
We continue the analysis of the entropy growth rate $h$ (\ref{eq:h_def}) of an typical random Markov matrix from the ensemble by decomposing it into an average term $h_{\rm ave}$ and an oscillating term $h_{\rm osc}$ reading
$$
  h = h_{\rm ave} + h_{\rm osc}\>,\qquad
  h_{\rm ave} = \frac{1}{d} \sum_{i=1}^d U_i\>,\quad
  h_{\rm osc} = \sum_{i=1}^d \left(\pi_i - \frac{1}{d}\right)(U_i - h_{\rm ave})\>,
$$
where we introduce auxiliary variables
$$
  U_i = - \sum_{j=1}^d M_{i,j} \log M_{i,j} \>.
$$
In the asymptotic limit $d\to\infty$ the variables $U_i$ have according to the CLT a Gaussian distribution with the first two central moments reading
\beqas
  \ave{U_i} &\sim& \log ( e^{\gamma -1} d)\>, \\
  \sigma_U^2 &\sim&  \left[1 + (\gamma-4) \gamma + \pi^2/3 + 
                      (2 \gamma - 4+\log d)\log d\right]\frac{1}{d} = O(d^{-1})\>,
\eeqas
with $\gamma$ being the Euler constant. The average $\ave{\cdot}$ in the expressions above is taken w.r.t. the asymptotic distribution $P_{\rm asym}(M)$ (\ref{eq:distr_mat_asym}). It is easy to see that the average term converges with increasing $d$ to $\ave{U_i}$ as
$$
  h_{\rm ave} = \ave{U_i} + O(d^{-\frac{1}{2}}) \>,
$$
where the last term on the r.h.s. denotes the statistical deviation. The oscillating term $h_{\rm osc}$ can be treated as a scalar product of vectors $(\pi_i - 1/d)_{i=1}^d$ and $(U_i-h_{\rm osc})_{i=1}^d$ and by applying the Schwarz-Cauchy inequality it be bounded from above:
$$
  h_{\rm osc}^2 \le \sum_{i=1}^d \left(\pi_i - \frac{1}{d}\right)^2  
  \left(U_i -h_{\rm ave}\right)^2 = O(d^{-2}) \>.
$$
The last term on r.h.s. denotes again the statistical deviation obtained by taking into account the statistical deviations $\ave{(\pi_i - 1/d)^2}=O(d^{-3})$ and $\ave{(U_i -h_{\rm ave})^2}=O(d^{-1})$. We conclude from this analysis that $h$ behaves over the ensemble, in the leading order in $d$, approximately as the Gaussian variable $h_{\rm ave}$ or $U_i$. This is supported by the numerical results in figure \ref{pic:entropy}.a, where we show that the cumulative distribution of $(h-\ave{h})/\sigma_h$ converges to the cumulative normal distribution $G(x)$ (\ref{eq:G_def}). By averaging $h$ over the ensemble the fluctuations from $\ave{U_i}$ are eliminated and we obtain the asymptotic equivalence
$$
  \ave{h} \sim  \log ( e^{\gamma -1} d) \quad \textrm{as} \quad d\to\infty \>.
$$
The later agrees very well with the numerically obtained $\ave{h}$ as we can see in figure \ref{pic:entropy}.b. This result supports the rule of thumb saying that in order to describe a dynamical system with an entropy rate $h$ accurately via a Markov process we need to a Markov matrix of the size $d\sim e^h$. 
The standard deviations $\sigma_h$ scales as $O(d^{-1})$, as seen in figure \ref{pic:entropy}.c, and this agrees with the expected statistical deviation of the form $\sqrt{(\sigma_U^2+h_{\rm osc}^2)/d}$ in our asymptotic approximation.

\section{The correlation decay induced by random Markov matrices}

\begin{figure}[!htb]
\centering
\includegraphics[width=6cm]{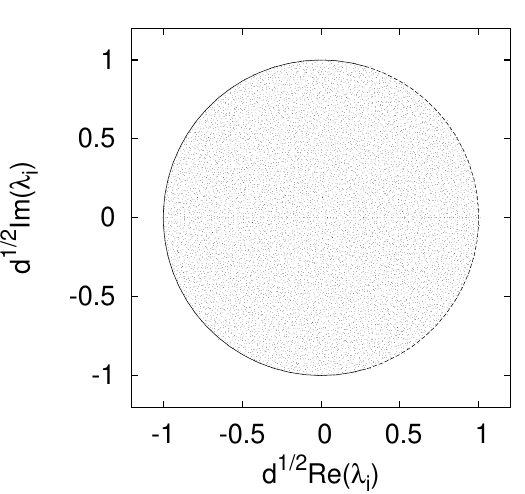}%
\includegraphics[width=8cm]{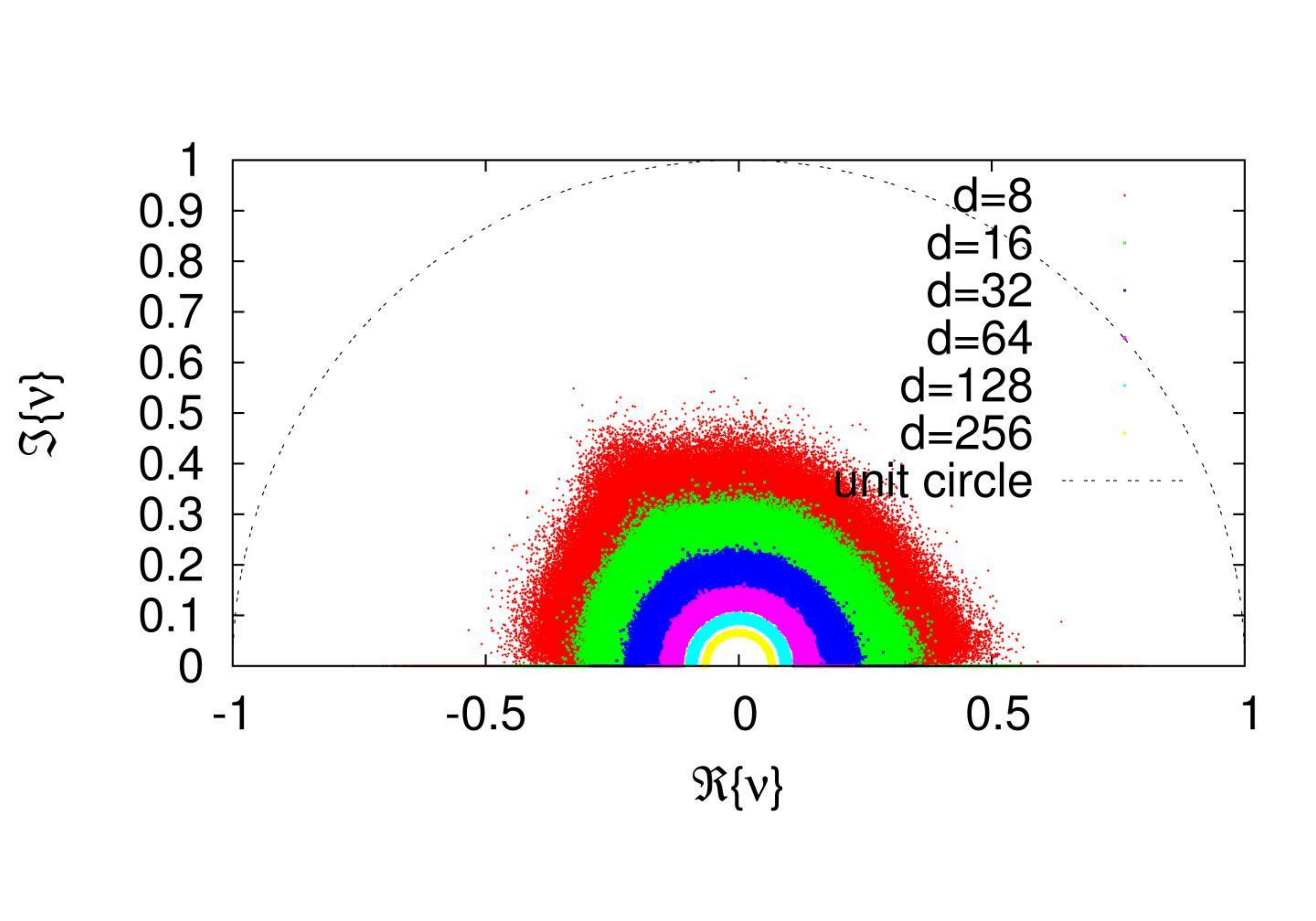}
\hbox to15cm{\hfil(a)\hfil\hfil(b)\hfil}
\caption{The spectrum of a random Markov matrix of dimension $d=10^4$ (a) without the eigenvalue $1$, which corresponds to the stationary distribution, and the second largest eigenvalues $\nu\in\bC$ (b) of approximately $N=10^6$ random Markov matrices calculated for different dimensions $d$.}
\label{pic:ruelle_points}
\end{figure}

\begin{figure}[!htb]
\centering 
\includegraphics[width=7cm]{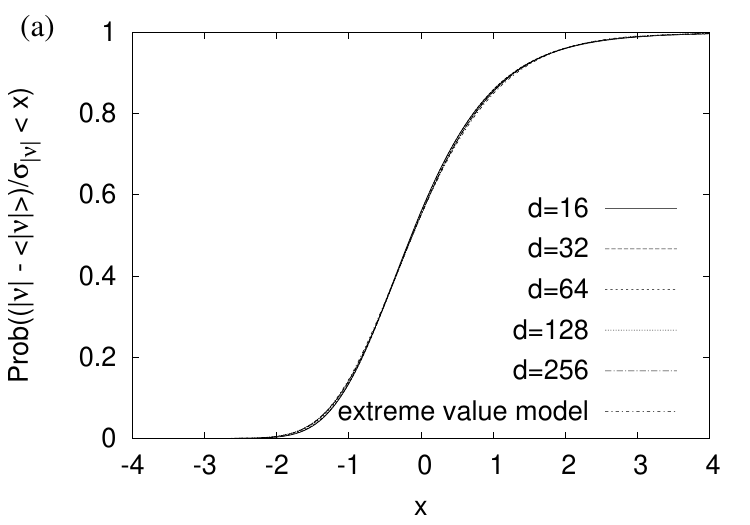}\\
\includegraphics[width=7cm]{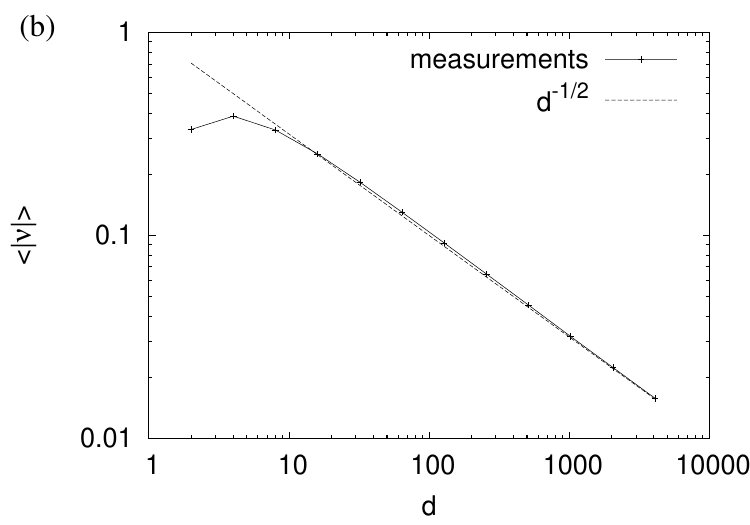}%
\includegraphics[width=7cm]{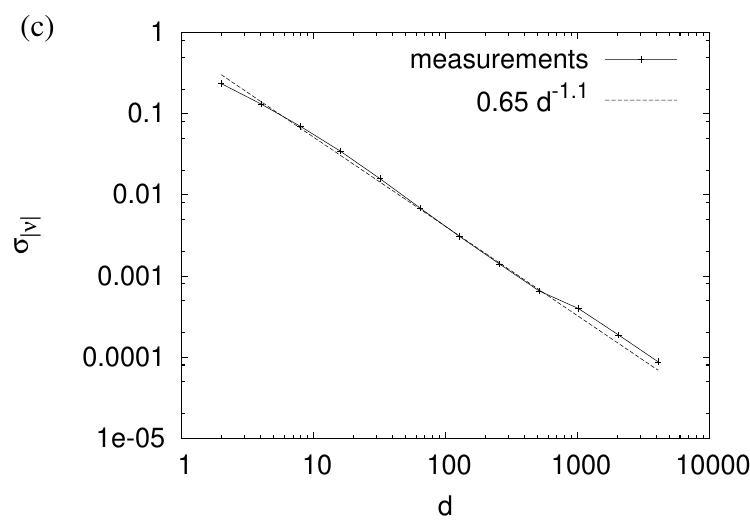}
\caption{The cumulative distribution of the magnitude of the second largest eigenvalue $|\nu|$ of random Markov matrices at different dimensions $d$ (a) and the average $\ave{|\nu|}$ (b) and the standard deviation of the amplitude $\sigma_{|\nu|}$ (b) as a function of the dimension $d$ calculated from approximately $10^6$ and $10^3$ random Markov matrices from the ensemble for $d<512$ and $d\ge 512$, respectively.}
\label{pic:ruelle_distr}
\end{figure}

\begin{figure}[!htb]
\centering
\includegraphics[width=7cm]{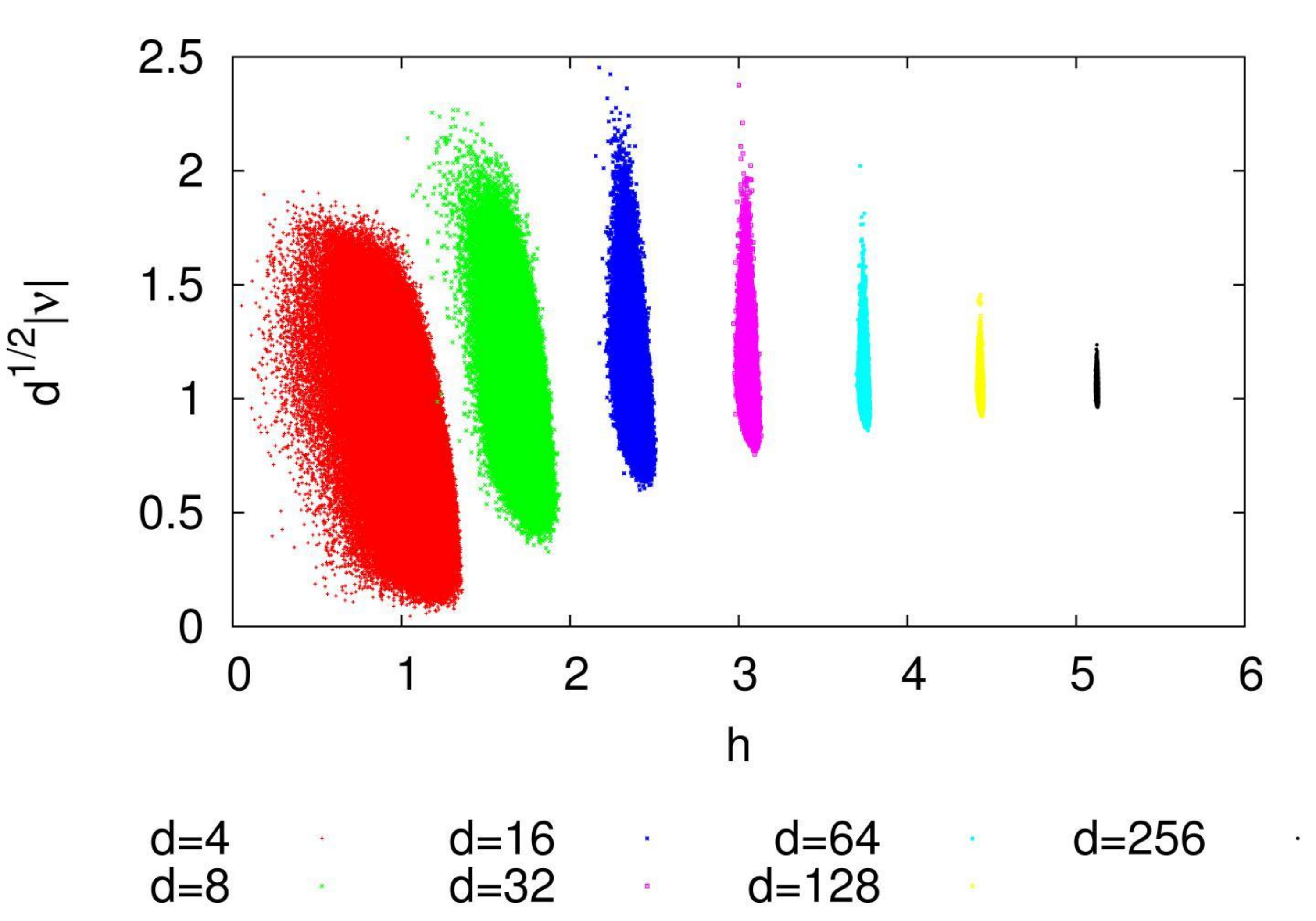}
\caption{The amplitude of the second largest eigenvalue $|\nu|$ and corresponding entropy $h$ calculated for some number of random Markov matrices samples uniformly with respect to the measure $\dd P(M)$ (\ref{eq:distr_mat}) at different dimensions $d$. For details on statistics see the caption under figure \ref{pic:ruelle_distr}.}
\label{pic:ent_ruelle}
\end{figure}

A state of a Markov chain defined by a Markov matrix $M=[M_{i,j}]_{i,j=1}^d\in\bR_+^{d\times d}$ is described by a probability distribution 
$$
  p=(p_i)_{i=1}^d \in \bR_+^d\>,\qquad 
  \sum_{i=1}^d p_i  = 1 \>,
$$
over given set of states $\{s_i\}_{i=1}^d$. Some initial probability distribution $p \in \bR_+^d $ is evolved in time to $p(t)$ by the Markov matrix in the following way 
$$
  {p(t)}^T = p^T M^t\>,
$$
where $t\in \bN_0$ denotes the discrete time. We find that a Markov chain generated by a typical random Markov matrix $M$ is mixing and consequently ergodic \cite{meyn:book:95}. We assume that the measure of Markov matrices in the ensemble corresponding to a non-mixing Markov chain is zero. \\
The discrete analogue of the time correlation function $C_{f,g}(t)$ between two real observables $f=(f_i\in \bR)_{i=1}^d$ and $g=(g_i\in \bR)_{i=1}^d$ is defined as
$$
  C_{f,g}(t) = \ave{f_i, (g(t))_i}_i - \ave{f_i}_i \ave{g_i}_i\>,
$$
where we introduce in time propagated observable ${g(t)}^T = g^T M^t$ and averaging over the stationary distribution $\ave{u_i}_i = \sum_i \pi_i u_i$. The second largest eigenvalue (called also the sub-dominant eigenvalue) of the Markov matrix $\nu \in \bC$ determines the decay of correlation between almost all pairs of observables $(f,g)$ following the formula 
$$
  |C_{f,g}(t)| = O(|\nu|^t)= O(e^{-t/\tau_{\rm c}})\quad\textrm{ as }\quad 
  t\to\infty\>,
$$
with $\tau_{\rm c} = - \log|\nu|$ called the correlation decay time. It is important to notice that spectrum $\Lambda=\{\lambda: \det(M-\lambda \id)=0\}$ of a Markov matrix $M$ has the symmetry
$$
  \Lambda^* = \Lambda\>,
$$
where $(\cdot)^*$ represents the complex conjugation. The symmetry can be noticed in figure (\ref{pic:ruelle_points}, where we show a spectrum of a typical random Markov matrix in the complex plane. In the limit of large dimensions $d\gg 1$ the eigenvalues are distributed symmetrically around the origin with the constant distribution of the square absolute value of the form
$$
  {\rm Prob}(x \le \|\lambda\|^2 < x+\dd x)  \approx 
  O(d^{-\frac{1}{2}})\,\dd x \>.
$$
This feature goes along the lines of the Girko's circular law \cite{girko:84} and its generalizations \cite{tao:07}, but this particular case is not yet proved to the best of our knowledge. Here we are mainly interested in the second largest eigenvalues $\nu\in\bC$ of the random Markov matrices. These are depicted for $N=10^6$ matrices sampled uniformly across the ensemble in figure \ref{pic:ruelle_points}.b for several different dimensions $d$.
For large $d$ the values of $\nu$ are distributed radially symmetrically around the origin with the average radius and dispersion decreasing with increasing $d$. Further we examine the distribution of the magnitude of the second largest eigenvalue $|\nu|$ denoted by $P_{|\nu|}$ and its first two central moments: average magnitude $\ave{|\nu|}$ and standard deviation $\sigma_{|\nu|}$. The cumulative distribution of the rescaled magnitude
$$
  \xi = \frac{|\nu| - \ave{\nu}}{\sigma_{|\nu|}}
$$
is depicted in figure \ref{pic:ruelle_distr}.a. From the figure we conclude that the distribution of $\xi$ is basically independent of dimension $d$ for large $d$ and we find that it agrees well with the extreme value statistics of type 1 (Gumbel) \cite{kotz:02}. Let us assume $x_i$ are i.i.d. standard Gaussian variables. Then the maximal value of $d$ variables $y = \max \{x_i\}_{i=1}^d$ is distributed according to the cumulative distribution
$$
  P_{\rm max}(y,d) = \left[G(y)\right]^d \>,
$$
where we use cumulative Gaussian distribution $G(y)$ (\ref{eq:G_def}). It is known that under simple linear transformation of the variable $y$, which depends on $d$, the distribution of transformed $y$ converges in the limit $d\to\infty$ to the Gumbel or double exponential distribution. To avoid certain slow convergence problems towards the limiting distribution outlined in \cite{resnik:87} we compare in figure \ref{pic:ruelle_distr}.a numerically obtained distribution directly with $P_{\rm max}((y-\mean{y})/\sigma_y,d)$ for several large enough $d$, where $\mean{y}$ and $\sigma_y$ are the average maximal value and its standard deviation, respectively. We find a very good agreement suggesting that the eigevalues of the Markov matrix behave as i.i.d random complex variables inside some disk in the complex plane with the radius $O(d^{-\frac{1}{2}})$. The first two central moments of numerical results as function of $d$ are shown in figures \ref{pic:ruelle_distr}.b and \ref{pic:ruelle_distr}.c.
The average magnitude of the second largest eigenvalue $\ave{|\nu|}$ fits very well to the  empirically found asymptotic formula
\beq
  \ave{|\nu|}\sim C_0\, d^{-\frac{1}{2}}\>,\qquad d\to\infty\>, 
  \label{eq:nu_asymp}
\eeq
where $C_0 \approx 1$. This result supports the conjecture that the Girko's circular law is valid for the random Markov matrices stated in \cite{chafai:07} and it can be better understood by the asymptotics of the upper bound for $|\nu|$ obtaining in following. By taking into account that all left-hand eigenvector except the stationary distribution $\pi$ are perpendicular to the vector $\one$, we can upper bound of the second largest eigenvalue $|\nu|$ as
$$
 |\nu|^2 \le \limsup_{x\in {\cal S}} \| x^\dag M \|_2^2\>,
 \qquad   
 {\cal S} =\{x \in \bC^d: \|x\|_2 = 1 \land x\perp \one\}\>.
$$
Here we write the expression $\| x^\dag M \|_2^2  = x^\dag N x$ using the matrix $N = M M^T$. In the asymptotic limit the matrix $N$ takes the form
$$
   N_{i,j} 
  =\sum_{k=1}^d M_{i,k} 
  M_{j,k} \sim \frac{1}{d}\delta_{i,j} + \one^T \one + O(d^{-2})\>.
$$
where the last term denotes the statistical error of the expression. From here we immediately obtain asymptotic expression of the upper bound
$$
  \limsup_{x\in {\cal S}} \| M x \|_2^2 \sim d^{-\frac{1}{2}} \>.
$$
This means that the second largest eigenvalue in a typical random Markov matrix is bounded by below $d^{-\frac{1}{2}}$ in the limit $d\to\infty$. This is true also in the average over the ensemble yielding the relation $\ave{\nu}\le d^{-\frac{1}{2}}$ and setting the value of the constant $C_0=1$. The asymptotic behaviour of the standard deviation $\sigma_{|\nu|}$ is not as clear as in the case of the average value $\ave{|\nu|}$. The numerical results suggest the power law decay
$$
  \sigma_{|\nu|} \sim C_1\, d^{-\alpha}\>,\qquad \alpha \approx 1.1\>.
$$
In the Markov approximations of a dynamical systems with the mixing property \cite{ott}, as described in the introduction, it interesting to know about the correlation decay-time $\tau_{\rm c}$ and entropy growth rate $h$ and their dependence on the cardinality of the state space. In the random Markov matrices form the ensemble we find that the average average correlation decay-time $\mean{t} :=-1/\log\ave{|\nu|}$ and average entropy growth rate $\ave{h}$ obey the asymptotics
$$
 \ave{h} \mean{\tau}_{\rm c}  \sim \frac{1}{2}\quad {\rm as} \quad d\to\infty \>.
$$
In dynamical systems there are strong indications that the correlation decay time and the entropy growth rate,  given as the sum of positive Lyapunov exponents, are correlated, but this connection is not well understood, yet, see e.g. \cite{ruelle, eckmann}. We address this question in the random Markov matrices and calculate pairs $(|\nu|, h)$ corresponding to the Markov matrices sampled uniformly over the ensemble. The result is graphically depicted in figure \ref{pic:ent_ruelle}, where one can notice that in particular at small dimensions $d$ there is clearly some correlation between the amplitude of the second largest eigenvalue $|\nu$ and entropy growth rate $h$ of random Markov matrices.

\begin{figure}[!htb]
\centering
\includegraphics[width=7cm]{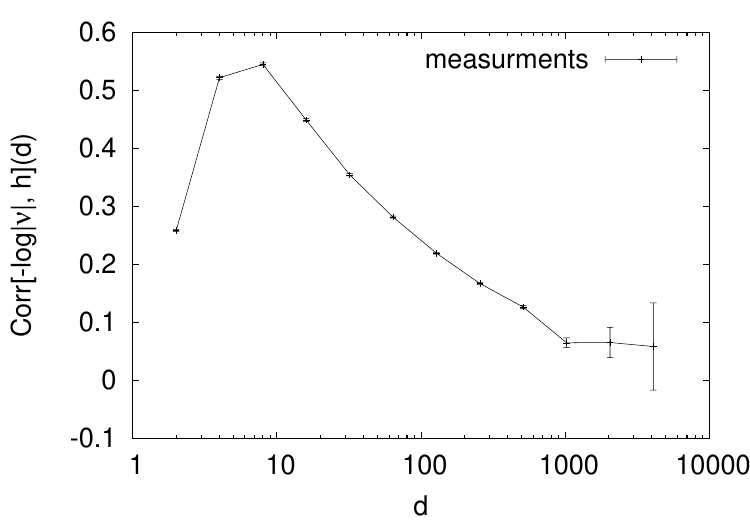}
\caption{The normalized correlation ${\rm Corr}(d)/{\rm Corr}(0)$ between the reciprocal correlation time $\tau^{-1}_{\rm c}=-\log|\nu|$ and the entropy growth rate $h$ in the ensemble of random Markov matrices as a function dimension $d$. For details on statistics see the caption under figure \ref{pic:ruelle_distr}.}
\label{pic:corr_mu_h}
\end{figure}

The latter is tested by calculating the statistical correlation between the reciprocal correlation decay time $\tau_{\rm c}^{-1}=-\log|\nu|$ and the entropy growth rate $h$ over the ensemble of random Markov matrices and is given by
$$
  {\rm Corr}(d)
  = \frac{\ave{ \left(\tau_{\rm c}^{-1}  - \ave{\tau_{\rm c}^{-1}}\right)
                \left(h -\ave{h}\right)}}
    {\sqrt{\ave{\left(\tau_{\rm c}^{-1} - \ave{\tau_{\rm c}^{-1}}\right)^2}
           \ave{\left(h -\ave{h}\right)^2}}} \>.
$$
The correlation ${\rm Corr}(d)$ as a function of dimension $d$ is presented in figure \ref{pic:corr_mu_h} and we see that it slowly decreases with increasing $d$. Currently it is not possible to determine the dependence of the correlation ${\rm Corr}(d)$  on dimension more precisely.

\section{Conclusions}

We define the ensemble of random Markov matrices, present its basic properties and point out few of the potential physical, technical and mathematical applications. We analyse the statistical properties of the stationary distribution $\pi=(\pi_i)_{i=1}^d$ corresponding to a typical element of the ensemble, and study the distribution of the entropy growth rate $h$ over the ensemble, where we find a good agreement with analytical predictions stating that $\pi_i$ is a Gaussian variable and $h$ is asymptotically equal to $\log(e^{\gamma-1} d)$ in the limit of large dimensions $d\to\infty$.
Further we analyse the second largest eigenvalue $\nu$ of the Markov matrices, which is connected to the correlation decay in the Markov chains. We show numerically and provide a heuristic proof that in average over the ensemble the second largest eigenvalue decreases with increasing dimension $d$ as $|\nu|\sim d^{-\frac{1}{2}}$ .
Additionally we calculate the correlation between the correlation decay rate and the entropy growth rate and find that it decreases with increasing dimension of the Markov matrices.\\
We believe that the current results enrich the understanding of Markov processes in the limit of large state spaces and all applications which can be described by Markov processes.

\section*{Acknowledgements}
The author would like to thank Karol \.Zyczkowski and Toma\v z Prosen for very interesting discussions and useful comments. The financial support by Slovenian Research Agency, grant Z1-0875-1554-08, is gratefully acknowledged.

\end{document}